\documentclass{amsart}
\usepackage{graphicx}
\usepackage{amssymb}

\newtheorem{thm}{Theorem}[section]

\newtheorem{prop}[thm]{Proposition}
\theoremstyle{definition}

\theoremstyle{remark}
\newtheorem{rem}[thm]{Remark}

\begin{document}

\title[Finite time blowup]{Finite time blowup for parabolic systems in two dimensions}
\author{Connor Mooney}
\address{Department of Mathematics, UT Austin, Austin, TX 78712}
\email{\tt cmooney@math.utexas.edu}


\begin{abstract}
We construct examples of finite time singularity from smooth data for linear uniformly parabolic systems in the plane. We obtain similar
examples for quasilinear systems with coefficients that depend only on the solution.
\end{abstract}
\maketitle

\section{Introduction}
We consider regularity for weak solutions to the linear parabolic system
\begin{equation}\label{ParabolicSystem}
{\bf u}_t = \text{div}(a(x,t)D{\bf u}).
\end{equation}
Here ${\bf u} : \mathbb{R}^n \times (-\infty, 0) \rightarrow \mathbb{R}^m$, and $a = [a^{ij}_{\alpha\beta}(x,t)]_{\alpha,\,\beta \leq m}^{i,\,j \leq n}$
are bounded measurable coefficients satisfying the uniform ellipticity condition
\begin{equation}\label{Ellipticity}
\lambda |p|^2 \leq a^{ij}_{\alpha\beta}(x,t)p^{\alpha}_ip^{\beta}_j \leq \Lambda |p|^2
\end{equation}
for some positive constants $\lambda,\,\Lambda$, and for all $p \in M^{m \times n}$ and all $(x,t)$. 
By a weak solution we mean a map ${\bf u} \in L^2_{loc}(\mathbb{R}^n \times (-\infty,0))$ with 
$D{\bf u} \in L^2_{loc}(\mathbb{R}^n \times (-\infty,0))$ that solves (\ref{ParabolicSystem}) in the sense of distributions. 
In coordinates one writes ${\bf u} = (u^1,...,u^m)$, and the system (\ref{ParabolicSystem}) is $u^{\alpha}_t = \partial_i( a^{ij}_{\alpha\beta}(x,t)u^{\beta}_j).$

Regularity results for (\ref{ParabolicSystem}) are important for the study of gradient flows in the calculus of variations. The gradient flow ${\bf v}$ of 
a functional with a smooth, uniformly convex integrand depending only on the gradient solves the system
\begin{equation}\label{GradientSystem}
 {\bf v}_t = \text{div}(B(D{\bf v})),
\end{equation}
where $B$ is a smooth uniformly monotone operator. The classical approach to regularity is to 
differentiate (\ref{GradientSystem}) and treat the problem as a linear system for the derivatives of ${\bf v}$ with bounded measurable coefficients.

Morrey \cite{Mo} showed that {\it stationary} solutions to (\ref{ParabolicSystem}) are continuous in the case $n = 2$. 
This follows from a higher-integrability result for the gradient. Solutions to (\ref{ParabolicSystem}) are also continuous
in the scalar case $m = 1$ by classical results of De Giorgi  \cite{DG1} and Nash \cite{Na}. As a consequence, solutions to (\ref{GradientSystem}) are smooth in these cases.
Solutions to (\ref{ParabolicSystem}) can be discontinuous in the case $n = m \geq 3$, by well-known examples of De Giorgi \cite{DG2} and Giusti-Miranda \cite{GM}. 

Ne\v{c}as and \v{S}ver\'{a}k \cite{NS} showed that time-dependent solutions to (\ref{GradientSystem}) are also smooth in the case $n = 2$.
However, in contrast with the scalar case and the planar elliptic case, the argument does not rely on continuity of solutions to the linearized problem. In fact, the question
of continuity of solutions to (\ref{ParabolicSystem}) in the case $n = 2$ remained open (stated e.g. in \cite{SJ} and \cite{JS}).
The purpose of this paper is to answer this question with a counterexample to regularity. Our main theorem is:

\begin{thm}\label{main}
 There exist a map
 $${\bf u} : \mathbb{R}^2 \times (-\infty, 0] \rightarrow \mathbb{R}^2$$
 that is smooth for $t < 0$ and Lipschitz up to $t = 0$ away from $(0,0)$, and a bounded matrix field
 $$a : \mathbb{R}^2 \times (-\infty, 0] \rightarrow \text{Sym}_{M^{2 \times 2} \times M^{2 \times 2}}$$
 satisfying (\ref{Ellipticity}), that is smooth for $t < 0$ and discontinuous at $(0,0)$, 
 such that ${\bf u}$ solves (\ref{ParabolicSystem}) in $\mathbb{R}^2 \times (-\infty, 0)$ with coefficients $a(x,t)$, and ${\bf u}(\cdot, 0)$ is discontinuous.
\end{thm}

\begin{rem}\label{unboundedmain}
 The example ${\bf u}$ in Theorem \ref{main} can in fact blow up in $L^{\infty}$.
\end{rem}

\begin{rem}
 One can extend to times $t \geq 0$ by e.g. keeping $a(x,t) = a(x,0)$ for $t > 0$, and solving the system with the initial data
 ${\bf u}(\cdot , 0)$. In this way one obtains a global (in space and time) weak solution that develops an interior discontinuity at $(0,0)$ which instantly
 disappears.
\end{rem}

\begin{rem}
For the system (\ref{ParabolicSystem}) there is a higher-integrability estimate for the spatial gradient in parabolic cylinders (see e.g. \cite{C}). In the case $n = 2$ 
this estimate implies that solutions are continuous in space at almost every time (which is not true when $n \geq 3$), but it does not rule out singularity formation.
\end{rem}

As a result of Theorem \ref{main}, one cannot rely on a continuity result at the linear level to prove regularity for (\ref{GradientSystem}) in the plane. 
One might instead hope to use that the derivatives of gradient flows solve quasilinear systems with the special structure
\begin{equation}\label{QuasilinearSystem}
 {\bf u}_t = \text{div}(a({\bf u})D{\bf u}),
\end{equation}
where $a^{ij}_{\alpha\beta}$ are smooth functions on $\mathbb{R}^m$ satisfying (\ref{Ellipticity}).
Our second result is an example of finite-time discontinuity from smooth data for the system (\ref{QuasilinearSystem}) in the case $n = 2, \, m = 4$:

\begin{thm}\label{nonlinearmain}
 There exist a map
 $${\bf u} : \mathbb{R}^2 \times (-\infty, 0] \rightarrow \mathbb{R}^4$$
 that is smooth for $t < 0$ and Lipschitz up to $t = 0$ away from $(0,0)$, and a smooth, bounded matrix field
 $$a : \mathbb{R}^4 \rightarrow \text{Sym}_{M^{4 \times 2} \times M^{4 \times 2}}$$
 satisfying (\ref{Ellipticity}),
 such that ${\bf u}$ solves (\ref{QuasilinearSystem}) in $\mathbb{R}^2 \times (-\infty, 0)$ with coefficients $a({\bf u})$, and ${\bf u}(\cdot, 0)$ is discontinuous.
\end{thm}

\begin{rem}
 The coefficients of the Giusti-Miranda example \cite{GM} can be written as smooth functions of ${\bf u}$, giving a discontinuous example in the case $n \geq 3$. 
\end{rem}

\begin{rem}
 It would be interesting to construct an example of finite time discontinuity from smooth data for (\ref{QuasilinearSystem}) in the case $n = m = 2$.
\end{rem}

Our examples show that parabolic systems in the plane behave differently than elliptic systems. They also show that the classical approach to proving regularity
for (\ref{GradientSystem}) in two dimensions fails. In \cite{NS} the authors instead prove a higher-integrability estimate for solutions of (\ref{ParabolicSystem}), and apply it to ${\bf v}_t$. 
One can then treat (\ref{GradientSystem}) as an elliptic system for each fixed time. 
Similar ideas were used to show the continuity of solutions to (\ref{ParabolicSystem}) in two dimensions when the coefficients are Lipschitz in space or in time
(see \cite{JS}).

The stationary examples of De Giorgi and Giusti-Miranda are discontinuous on the cylindrical set $\{x = 0\}$. Examples of finite time discontinuity from smooth data for (\ref{ParabolicSystem})
were constructed in the case $n = m \geq 3$ by Star\'{a}, John and Mal\'{y} in \cite{SJM}, and refined
by Star\'{a} and John in \cite{SJ}. In these examples, the data and coefficients are a small perturbation from those of the De Giorgi example.

The data in our examples are also a perturbation of the De Giorgi example, but due to low-dimensionality we need to take a different approach to constructing the coefficients,
and also to make a more careful perturbation. To prove Theorem \ref{main} we search for a solution of the form ${\bf u} = {\bf U}(x/\sqrt{-t})$. This reduces
the problem to finding a nontrivial global, bounded solution to an elliptic system. Our approach is to construct a pair of functions that solve the analogous scalar equation away from an annulus,
where the error in the equation is small. This pair defines a map that solves a decoupled system away from the annulus. 
We then couple the equations so that the system is solved globally.

\begin{rem}\label{LiouvilleRemark}
 An important feature of our example is that $|{\bf U}|$ is not radially increasing, unlike in the higher-dimensional examples.
 In fact, such examples do not exist in the plane. In Section \ref{Appendix} we prove a Liouville theorem in two dimensions for self-similar solutions 
 with radially increasing modulus (see Theorem \ref{Liouville}).
\end{rem}

Our remaining examples are modifications of the construction described above.
To obtain a solution to (\ref{ParabolicSystem}) with $L^{\infty}$ blowup we instead search for solutions invariant under rescalings that fix $-\epsilon$-homogeneous maps.

Because $|{\bf u}|$ is not radially increasing in our first example (which is guaranteed by the Liouville theorem mentioned in Remark \ref{LiouvilleRemark}),
we can not write the coefficients as functions of ${\bf u}$ (see Remark \ref{UDependence}). 
To prove Theorem \ref{nonlinearmain} we go to higher codimension. 
We take a solution $\tilde{\bf u}$ to (\ref{ParabolicSystem}) that is similar to ${\bf u}$, such that the map $|x| \rightarrow (|{\bf u}|,\,|\tilde{\bf u}|)$ is injective.
The pair $({\bf u}, \, \tilde{\bf u})$ solves a uniformly parabolic system in the case $n = 2,\, m = 4$, and we can write the coefficients as smooth functions of $({\bf u},\,\tilde{\bf u})$.

The paper is organized as follows. In Section \ref{Reduction} we reduce Theorem \ref{main} to finding a global, bounded solution ${\bf U}$ to an elliptic system by searching for solutions 
that are invariant under parabolic scaling.
In Section \ref{BuildingBlock} we construct a function that solves the analogous elliptic equation away from an annulus.
Using this function we define ${\bf U}$ and diagonal coefficients so that ${\bf U}$ solves the desired (decoupled) system away from the annulus.
In Section \ref{Coupling} we construct off-diagonal coefficients that couple the equations so that ${\bf U}$ solves the system globally, and we verify that the resulting matrix field is uniformly elliptic. 
This completes the proof of Theorem \ref{main}.
In Section \ref{Unbounded} we modify this construction to obtain an example with $L^{\infty}$ blowup.
In Section \ref{QuasilinearStructure} we prove Theorem \ref{nonlinearmain}.
Finally, in Section \ref{Appendix} we prove a Liouville theorem indicating why $|{\bf U}|$ can not be radially increasing
in two dimensions.

\section{Reduction}\label{Reduction}
We first reduce the problem to
finding a global bounded solution to an elliptic system by searching for solutions that are invariant under the parabolic scaling 
$(x,t) \rightarrow (\lambda x, \lambda^2 t)$.

\begin{prop}\label{EllipticReduction}
Assume that ${\bf U} : \mathbb{R}^n \rightarrow \mathbb{R}^m$ is a non-constant, bounded, smooth solution to the system
\begin{equation}\label{StationarySystem}
 \text{div}(A(x)D{\bf U}) = \frac{1}{2}D{\bf U} \cdot x,
\end{equation}
where $A = A^{ij}_{\alpha\beta}(x)$
are smooth, uniformly elliptic coefficients. If we take
$${\bf u}(x,t) := {\bf U}\left(\frac{x}{\sqrt{-t}}\right), \quad a(x,t) = A\left(\frac{x}{\sqrt{-t}}\right),$$ 
then ${\bf u}$ solves (\ref{ParabolicSystem}) on $\mathbb{R}^n \times (-\infty, 0)$ with the coefficients $a(x,t)$.

Furthermore, if ${\bf U}$ satisfies
\begin{equation}\label{BoundaryConditions}
|D{\bf U}(x)| = O(|x|^{-1}), \quad |D{\bf U} \cdot x| = O(|x|^{-2}),
\end{equation}
then ${\bf u}$ is smooth for $t < 0$ and Lipschitz up to $t = 0$ away from $(0,0)$, and is discontinuous at $(0,0)$.
\end{prop}

The proof is a straightforward computation.

\begin{rem}
To produce an example with $L^{\infty}$ blowup we instead search for solutions of the form $(-t)^{-\epsilon/2}{\bf U}(x/\sqrt{-t})$, where
${\bf U}$ satisfies estimates analogous to (\ref{BoundaryConditions}) at infinity (see Section \ref{Unbounded}).
\end{rem}

\begin{rem}\label{QuasilinearReduction}
Likewise, if ${\bf U}$ solves $\text{div}(A({\bf U})D{\bf U}) = \frac{1}{2}D{\bf U} \cdot x$ where $A$ are smooth uniformly elliptic coefficients on $\mathbb{R}^m$, then ${\bf u}(x,t) = {\bf U}\left(\frac{x}{\sqrt{-t}}\right)$
solves (\ref{QuasilinearSystem}) on $\mathbb{R}^n \times (-\infty, 0)$ with coefficients $a({\bf u}) = A({\bf U})$.
\end{rem}

\begin{rem}\label{ReductionLiouville}
 The problem of finding self-similar singular solutions to (\ref{ParabolicSystem}) thus boils down to proving or disproving a Liouville theorem 
 for the system (\ref{StationarySystem}).
 In Section \ref{Appendix} we verify the Liouville theorem in the case that $|{\bf U}|$ is radially increasing and $n = 2$.
\end{rem}


\section{Scalar Building Block}\label{BuildingBlock}
We now construct a smooth function $u : \mathbb{R}^2 \rightarrow \mathbb{R}$ and a smooth, uniformly elliptic matrix 
field $M: \mathbb{R}^2 \rightarrow \text{Sym}_{2 \times 2}$ such that $u$ solves
\begin{equation}\label{ScalarEquation}
\frac{1}{2}\nabla u \cdot x - \text{div}(M\nabla u) = 0
\end{equation}
away from an annulus, where the expression on the left side is small.

For $x$ in the plane, we denote $|x|$ by $r$ and the unit radial and tangential vectors $\nu$ and $\tau$ by
$$\nu = \frac{x}{r}, \quad \tau = \frac{x^{\perp}}{r}$$
away from the origin, where $x^{\perp}$ is the counterclockwise rotation of $x$ by $\frac{\pi}{2}$. Observe that
\begin{equation}\label{DivZero}
\text{div}\left(\frac{\nu}{r}\right) = \text{div}\left(\frac{\tau}{r}\right) = 0 
\end{equation}
away from the origin, since they are the gradients of harmonic functions.

Now let
$$u = \varphi(r)\cos(\theta)$$
and
$$M = f(r)\nu \otimes \nu + h(r) \tau \otimes \tau$$
for some $\varphi$ and positive bounded $f,\, h$ to be chosen. The left side of Equation (\ref{ScalarEquation}) can be written $$E(r)\cos\theta,$$ where 
\begin{equation}\label{EquationDeficit}
 E(r) := \frac{1}{2}r\varphi' + \frac{h\varphi}{r^2} - \frac{(r\varphi'f)'}{r}.
\end{equation}
This follows from a short computation using (\ref{DivZero}) and that
$$\nabla u = r\varphi'(r) \cos\theta \, \frac{\nu}{r} - \varphi(r) \sin\theta \, \frac{\tau}{r}.$$

\subsection{Definition of $\varphi$}
\begin{figure}
 \centering
    \includegraphics[scale=0.35]{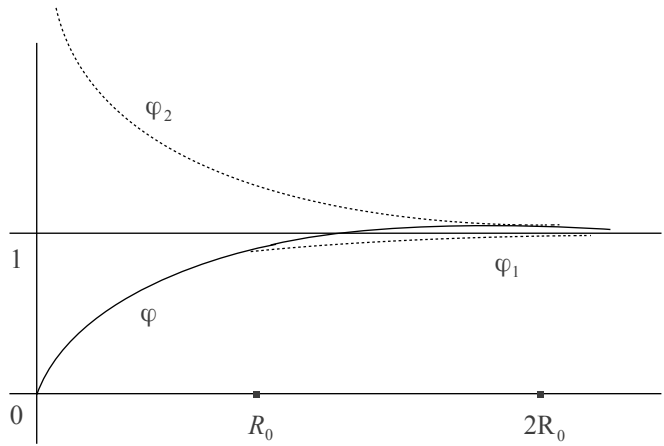}
 \caption{The function $\varphi$ smoothly connects $\varphi_1$ and $\varphi_2$ on $[R_0,\,2R_0]$, and satisfies
 the estimates $|\varphi'| < Cr^{-3},\,|\varphi''| < Cr^{-4}$.}
 \label{phi_construction}
\end{figure}

Define
$$\varphi_1 = \frac{r}{\sqrt{1 + r^2}}, \quad \varphi_2 = 1 + \frac{1}{2r^2}.$$
Let $\xi$ be a smooth, non-increasing function that is $1$ to the left of zero and $0$ to the right of one. For some $R_0$ large to be chosen let
$$\varphi(r) = \xi\left(\frac{r-R_0}{R_0}\right)\varphi_1 +  \left(1-\xi\left(\frac{r-R_0}{R_0}\right)\right)\varphi_2$$
(See Figure \ref{phi_construction}).

The following estimates are easy to verify:
\begin{equation}\label{phiEstimate}
 \varphi'(r) \leq Cr^{-3}, \quad \varphi''(r)  \leq Cr^{-4}.
\end{equation}
(Here and below $C$ denotes a universal constant independent of $R_0$).

\begin{rem}
 The motivation for our choice of $\varphi$ is as follows. We want $u$ to look $0$-homogeneous for $r$ large, so the angular derivatives dominate and one has 
 $\Delta u \sim -r^{-2}u$. Thus, solving the heat equation
 with initial data $u$ is compatible with ``squeezing'' by parabolic rescaling if $\varphi$ is decreasing at the rate $r \varphi' \sim -r^{-2}$.
 One can solve the equation $E(r) = 0$ where $\varphi' > 0$ by letting the coefficient $f$ grow large (see below), 
 but near the circle $\{\varphi' = 0\}$ the function $u$ can not solve the desired equation by the maximum principle.
 
\end{rem}

\subsection{Definition of $f$ and $h$}
\begin{figure}
 \centering
    \includegraphics[scale=0.35]{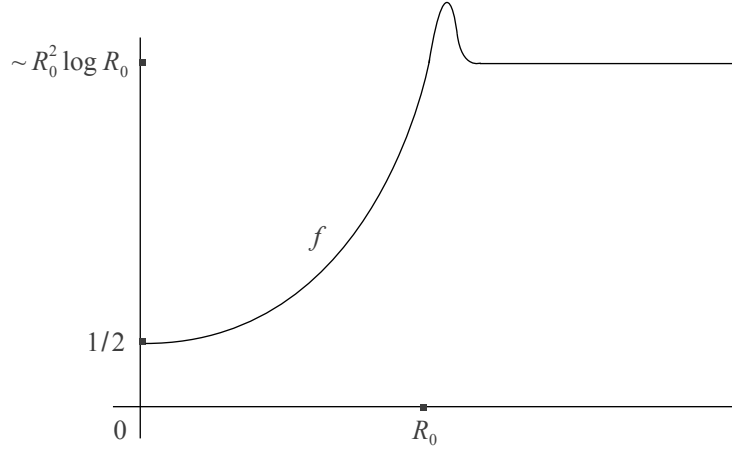}
 \caption{The function $f$ increases from $1/2$ to $\sim R_0^2\log R_0$ on $[0,R_0]$, then remains constant.}
 \label{f_construction}
\end{figure}

For $r < R_0$ we can solve the equation $E(r) = 0$
by keeping $h$ bounded and allowing $f$ to grow. Taking $h = 1/2$ for $r < R_0$ and solving $E(r) = 0$ for $f$ gives the function
\begin{align*}
f_0(r) &= \frac{(1 + r^2)^{3/2}}{2}\,\frac{1}{r} \int_0^r \frac{1 + 2s^2}{(1 + s^2)^{3/2}} \,ds \\
&= \frac{(1+r^2)^{3/2}}{r}\log((1+r^2)^{1/2} + r) - \frac{1}{2}(1 + r^2).
\end{align*}
It is straightforward to check that $f_0$ is strictly positive and locally bounded, and that the expansion of $f_0$ around $0$ has only even powers of $r$ (so its even reflection is smooth).
Furthermore, $f_0$ has the asymptotics
\begin{equation}\label{fEstimates}
 R^2 \log R \leq f_0(R) \leq 2R^2\log R, \quad R > R_0
\end{equation}
for $R_0$ sufficiently large. We take
$$f(r) := f_0(r)\xi(r - R_0) + (1-\xi(r-R_0))f_0(R_0)$$
(see Figure \ref{f_construction}).

\begin{figure}
 \centering
    \includegraphics[scale=0.35]{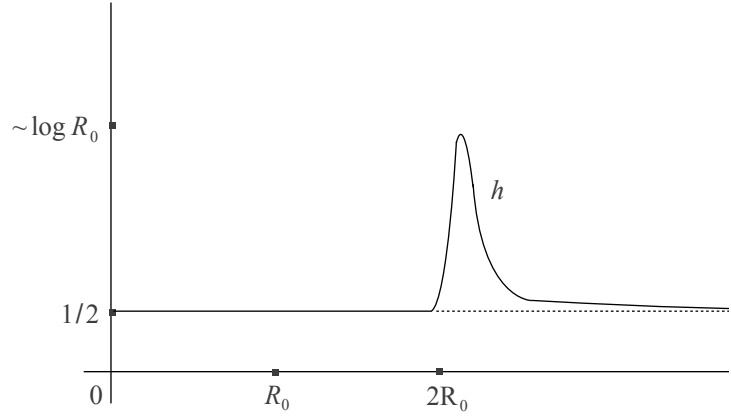}
 \caption{The function $h$ is close to $1/2$ most of the time, with a bump near $2R_0$ so the equation is solved for $r > 2R_0 + 1$.}
 \label{h_construction}
\end{figure}

Now define
$$h_0 := 1/2, \quad  \quad h_1 := \frac{1}{\varphi}\left(\frac{1}{2} + \frac{2f(R_0)}{r^2}\right).$$
One checks using the definition of $f$ and $\varphi$ that for $r > 2R_0$, one has $E(r) = 0$ by taking $h = h_1$. We define
$$h(r) =  \xi(r - 2R_0)h_0 + (1-\xi(r-2R_0))h_1$$
(see Figure \ref{h_construction}). Note that $h$ satisfies 
\begin{equation}\label{hEstimate}
 1/2 \leq h \leq C\log R_0.
\end{equation}

\begin{figure}
 \centering
    \includegraphics[scale=0.35]{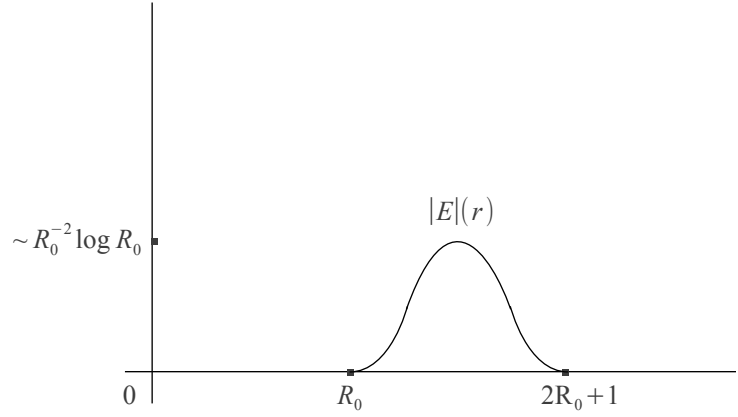}
 \caption{The error in the equation is supported in $[R_0,\, 2R_0 + 1]$ and is of order $R_0^{-2}\log R_0.$}
 \label{E_Deficit}
\end{figure}

With these choices of $f,\,h$, we have that
$$E(r) = 0, \quad r \in [R_0, \, 2R_0 + 1].$$
By the estimates (\ref{phiEstimate}), (\ref{fEstimates}) and (\ref{hEstimate}), in the remaining annulus we have
\begin{equation}\label{DeficitEstimate}
 |E(r)| \leq C\left(\frac{\log R_0}{r^2} + \frac{R_0^2\log R_0}{r^4}\right) \, \chi_{[R_0,\,2R_0 + 1]} < CR_0^{-2}\log R_0 \, \chi_{[R_0, \,2R_0 + 1]}
\end{equation}
(see Figure \ref{E_Deficit}).

Furthermore, one checks for $r < R_0$ that
$$M = \frac{1}{2} I + \beta(r) x \otimes x$$
where $\beta(|x|) = \frac{f(|x|) - 1/2}{|x|^2}$ is a smooth function on $B_{R_0}$. Thus,
$M$ is smooth, bounded and uniformly elliptic on $\mathbb{R}^2$ with eigenvalues between $\frac{1}{2}$ and $CR_0^2\log R_0$.

\subsection{Definition of ${\bf U}$}
We define the components of ${\bf U}$ by $u$ and a rotation of $u$:
$${\bf U} = (u^1,\,u^2) = (\varphi(r)\cos\theta, \, \varphi(r)\sin\theta) = \varphi(r)\nu.$$
Using the estimates (\ref{phiEstimate}) for $\varphi$ one verifies that
\begin{equation}\label{UEstimates}
 |D{\bf U}| = O(r^{-1}), \quad |D{\bf U} \cdot x| = O(r^{-2})
\end{equation}
as desired.

Furthermore, taking $B_{11} = B_{22} = M$ and $B_{12} = B_{21} = 0$, by construction and the rotation invariance of $M$ the map ${\bf U}$ solves the equation 
$$\frac{1}{2}D{\bf U} \cdot x - \text{div}(BD{\bf U}) = E(r)\nu.$$
In the next section we will perturb $B_{12}$ and  $B_{21}$ so that the system is solved globally and the coefficients remain uniformly elliptic.

\section{Coupling the Equations}\label{Coupling}
By the analysis above, if we take $A_{11} = A_{22} = M$ and $A_{12} = A_{21} = 0$, then
the map ${\bf U}$ solves the desired elliptic equation (\ref{StationarySystem}) away from the annulus $R_0 < r < 2 R_0 + 1$.
We now couple the equations in this region. We will use that $f(r)$ is large in the annulus to 
conclude that the resulting coefficient matrix $A$ is uniformly elliptic.

\begin{figure}
 \centering
    \includegraphics[scale=0.35]{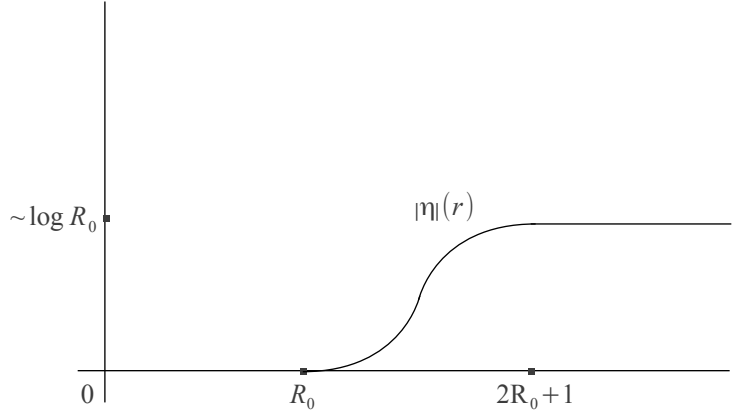}
 \caption{The corrector $\eta$ is zero to the left of $R_0$ and constant to the right of $2R_0+1$, with $|\eta|$ of order $\log R_0$.}
 \label{corrector}
\end{figure}

Since $u^2$ is a rotation of $u^1$ is natural to look for coupling coefficients that are rotations. Let $A_{12}$ be the ``corrector'' matrix field
$$A_{12} = \eta(r) \left(\begin{array}{cc}
 0 & 1 \\
 -1 & 0
 \end{array}\right).$$
One computes
$$\text{div}(A_{12} \nabla u^2) = \frac{\eta'\varphi}{r}\cos\theta.$$
Thus, to solve the equation (\ref{StationarySystem}) we need to take
$$\eta(r) := \int_0^r \frac{tE(t)}{\varphi(t)}\,dt.$$
With this choice of $\eta$, the desired equation 
$$\text{div}(A_{11}\nabla u^1 + A_{12}\nabla u^2) = \frac{1}{2}\nabla u^1 \cdot x$$ 
is solved, and by the estimate (\ref{DeficitEstimate}) we have
\begin{equation}\label{CorrectorEstimate}
 |\eta(r)| \leq C\log R_0 \, \chi_{\{r > R_0\}}
\end{equation}
(see Figure \ref{corrector}).

Finally, we define the remaining corrector $A_{21}$ by
$$A_{21} = -A_{12},$$
so that the equation holds in the second component.

In conclusion, we constructed a coefficient matrix $A$ and a map ${\bf U}$ solving the system (\ref{StationarySystem}). With respect to the coordinate system
$$(\nu,0),\, (\tau,0), \, (0,\nu), \,(0,\tau)$$
(where $(v,w)$ denotes the $2 \times 2$ matrix with first row $v$ and second row $w$) one writes
$$ A = 
\left(\begin{array}{cc|cc}
 f & 0 & 0 & \eta \\
 0 & h & -\eta & 0 \\
 \hline
 0 & -\eta & f & 0 \\
 \eta & 0 & 0 & h \\
\end{array}\right)(r).
$$

For $r < R_0$ one has $\eta = 0$ and the equations are decoupled. For $r > R_0$ large we examine the characteristic polynomial
$$P(\lambda) = \left[(\lambda -f)(\lambda - h) - \eta^2\right]^2.$$
Using the estimates (\ref{fEstimates}), (\ref{hEstimate}) and (\ref{CorrectorEstimate})
one sees that, for $\lambda \leq 0$, we have $P(\lambda) > 0$, verifying uniform ellipticity and completing the example:

\begin{proof}[{\bf Proof of Theorem \ref{main}}]
The map ${\bf U}$ and matrix field $A$ satisfy the hypotheses of Proposition \ref{EllipticReduction} by construction and estimate (\ref{UEstimates}).
\end{proof}

\begin{rem}\label{UDependence}
It is not hard to write $f$ as a smooth function of $\varphi$ and $h$ as a Lipschitz function of $\varphi$. However, on the circle $\{\varphi' = 0\}$, one computes that $E > 0$. (Indeed, the error
must be nonzero there by the maximum principle). It follows that 
$\eta$ is not a function of $\varphi$. In particular, the coefficients cannot be written as functions of ${\bf U}$. We overcome this in Section \ref{QuasilinearStructure} by going to higher codimension.
\end{rem}

\section{Unbounded Singularity}\label{Unbounded}
In this section we modify the construction from the previous section to produce an example with $L^{\infty}$ blowup at $(0,0)$. 
The construction follows the same lines, so we just sketch the key steps. For simplicity we use the same notation as above.

\vspace{2mm}

{\bf Reduction to Elliptic System.}
We search for solutions of the form
$${\bf u}(x,t) = \frac{1}{(-t)^{\epsilon/2}}{\bf U}\left(\frac{x}{\sqrt{-t}}\right)$$
for some $\epsilon > 0$, with coefficients
$$a(x,t) = A\left(\frac{x}{\sqrt{-t}}\right).$$
The idea is that this rescaling fixes $-\epsilon$-homogeneous functions rather than $0$-homogeneous functions.
This reduces the problem to finding a nontrivial smooth, global bounded solution ${\bf U}$ to the elliptic system
\begin{equation}\label{UnboundedStationarySystem}
\text{div}(AD{\bf U}) = \frac{1}{2}(D{\bf U} \cdot x + \epsilon {\bf U}),
\end{equation}
where $A(x)$ are smooth uniformly elliptic coefficients and ${\bf U}$ satisfies
\begin{equation}\label{UnboundedDesiredEstimates}
 |D{\bf U}| = O(|x|^{-1-\epsilon}), \quad |D{\bf U} \cdot x + \epsilon{\bf U}| = O(|x|^{-2-\epsilon}).
\end{equation}
One checks that if ${\bf U}$ satisfies these conditions, then ${\bf u}$ is smooth for $t < 0$ and Lipschitz up to $t = 0$ away from $(0,0)$ and
$\|{\bf u}(\cdot ,t)\|_{L^{\infty}(B_1)}$ blows up at the rate $(-t)^{-\epsilon/2}$.

\begin{rem}
 In fact, we will choose ${\bf U}$ to be asymptotically homogeneous of degree $-\epsilon$, so that ${\bf u}(\cdot ,0)$ is homogeneous of degree $-\epsilon$.
\end{rem}

\vspace{2mm}

{\bf Scalar Building Block.} We will again build ${\bf U}$ out of a scalar function $u$ that solves the elliptic equation
$$\text{div}(M \nabla u) = \frac{1}{2} (\nabla u \cdot x + \epsilon u)$$
away from an annulus. Take
$$u = \varphi(r)\cos\theta, \quad M = f(r) \nu \otimes \nu + h(r) \tau \otimes \tau.$$ 
In this case we have 
$$\frac{1}{2}(\nabla u \cdot x + \epsilon u) - \text{div}(M\nabla u) = E(r)\cos\theta$$ 
with
$$E(r) := \frac{1}{2}(r\varphi' + \epsilon \varphi) + \frac{h\varphi}{r^2} - \frac{(r\varphi'f)'}{r}.$$

\vspace{2mm}

{\bf Definition of $\varphi$.} 
We take $\varphi = \varphi_1$ (the same as above) for $r < R_0$ large, and for $r > 2R_0$ we define
$$\varphi(r) = \varphi_3(r) := r^{-\epsilon} + \frac{1}{2}r^{-\epsilon-2}.$$
Note that for $\epsilon = 0$ this reduces to what we have above. Take
$$\epsilon = \frac{1}{R_0^2\log R_0}.$$ 
Then in the interval $[R_0, \, 2R_0]$ one verifies
$$|\varphi_3'| < CR_0^{-3}, \quad |\varphi_3''| < CR_0^{-4}.$$
Furthermore, since $1- R_0^{-\epsilon} \leq C \epsilon \log R_0 \leq C R_0^{-2}$, we can take $\varphi$ to be a smooth
gluing of $\varphi_1$ to $\varphi_3$ in $[R_0,\,2R_0]$ so that same estimates as above hold in the corrector region:
\begin{equation}\label{UnboundedPhiEstimate}
 |\varphi'| < \frac{C}{R_0^3}, \quad |\varphi''| < \frac{C}{R_0^4} \quad \text{ for } R_0 \leq r \leq 2R_0.
\end{equation}

\vspace{2mm}

{\bf Construction of $f$ and $h$.}
Take $h = 1/2$ for $r < R_0$ and solve $E(r) = 0$ for a function $f_0$. Then $f_0(|x|)$ is positive and smooth for $|x| < R_0$ with the asymptotics
\begin{equation}\label{UnboundedfEstimate}
 f_0(R_0) \sim R_0^2 \log R_0 + \epsilon R_0^4 \sim R_0^2 \log R_0.
\end{equation}
(Here $\sim$ denotes equivalence up to multiplying by constants independent of $R_0$).
Define $f$ to be a gluing of $f_0$ to $f_0(R_0)$ between $R_0$ and $R_0 + 1$ as above.

We again choose $h$ so that $E(r) = 0$ for $r > 2R_0 + 1$. The error in $\{r > 2R_0\}$ is
$$E(r) = r^{-2-\epsilon}\left(-\frac{1}{2} + \left(1 + \frac{1}{2}r^{-2}\right)h - f(R_0)\left(\epsilon^2 + \frac{(2+\epsilon)^2}{2}r^{-2}\right)\right).$$
So we define $h$ in $\{r > 2R_0 + 1\}$ by
$$(1 + r^{-2}/2)h(r) = \frac{1}{2} + f(R_0)\left(\epsilon^2 + \frac{(2+\epsilon)^2}{2}r^{-2}\right),$$
and glue it to $1/2$ for $r < 2R_0$. This gives 
\begin{equation}\label{UnboundedhEstimate}
 \frac{1}{2} \leq h \leq C\log R_0,
\end{equation}
with $h$ asymptotically close to $1/2$ and with a bump of size $\log R_0$ near $2R_0$.

\vspace{2mm}

{\bf Definition of ${\bf U}$.} We again let
$${\bf U} = \varphi(r)\nu.$$
One checks using the definition of $\varphi$ that the derivatives of ${\bf U}$ satisfy the desired estimates (\ref{UnboundedDesiredEstimates}).
If we take $B_{11} = B_{22} = M$ and $B_{12} = B_{21} = 0$ then ${\bf U}$ solves
$$\frac{1}{2}(D{\bf U} \cdot x + \epsilon {\bf U}) - \text{div}(BD{\bf U}) = E(r)\nu,$$
and using the estimates (\ref{UnboundedPhiEstimate}), (\ref{UnboundedfEstimate}) and (\ref{UnboundedhEstimate}) we conclude that the error is estimated by
\begin{equation}\label{UnboundedErrorEstimate}
 |E(r)| \leq C\frac{\log R_0}{R_0^2} \, \chi_{[R_0, \, 2R_0+1]}.
\end{equation}

\vspace{2mm}

{\bf Coupling the equations.}
Let $A_{11} = M$ and again take
$$A_{12} = \eta(r) \left(\begin{array}{cc}
 0 & 1 \\
 -1 & 0
 \end{array}\right).$$
To solve the desired equation
$$\text{div}(A_{11}\nabla u^1 + A_{12} \nabla u^2) = \frac{1}{2}(\nabla u^1 \cdot x + \epsilon u^1)$$ 
we again need
$$\frac{\eta'\varphi}{r} = E(r).$$
Integrating and using (\ref{UnboundedErrorEstimate}) we obtain
\begin{equation}\label{UnboundedCorrectorEstimate}
 |\eta| \leq C\log R_0 \, \chi_{\{r > R_0\}}.
\end{equation}

Taking $A_{22} = M$ and $A_{21} = -A_{12}$ one verifies that the desired system (\ref{UnboundedStationarySystem}) is also solved in the second component.
Finally, the resulting matrix $A$ is smooth, and the estimates (\ref{UnboundedfEstimate}), (\ref{UnboundedhEstimate}) and (\ref{UnboundedCorrectorEstimate})
give that $A$ is positive, completing the example.

\begin{rem}
 In the above construction we see that $\|u(\cdot, t)\|_{L^{\infty}(B_1)}$ blows up at the rate 
 $(-t)^{-\frac{1}{2}R_0^{-2}(\log R_0)^{-1}}.$ A natural question is how quickly a solution to (\ref{ParabolicSystem}) in two dimensions can blow up in $L^{\infty}$
 from smooth data, i.e. how large one can take $\epsilon$.
\end{rem}

\begin{rem}
 We remark that our examples are smooth for $t < 0$. In \cite{SJ} the authors construct an example with finite time blowup in the case 
 $n = m \geq 3$ that is H\"{o}lder continuous, but not smooth, for $t < 0$.
\end{rem}

\section{An Example for Quasilinear Structure}\label{QuasilinearStructure}
In this section we construct a solution to the quasilinear problem (\ref{QuasilinearSystem}) that develops an interior discontinuity in finite time from smooth data.
We will construct a smooth, bounded map ${\bf W} : \mathbb{R}^2 \rightarrow \mathbb{R}^4$ and smooth matrix field $A({\bf W})$
satisfying the hypotheses in Remark \ref{QuasilinearReduction}, and the estimates (\ref{BoundaryConditions}).

\subsection{Construction of {\bf W}}
Let ${\bf U}$ be the map constructed in Section \ref{BuildingBlock}.
Recall that ${\bf U} = \varphi(r)\nu$ where $\varphi(r)$ smoothly connects $\varphi_1$ to $\varphi_2$ in the interval $[R_0,\,2R_0]$.
We let $\tilde{\bf U} = \tilde{\varphi}(r)\nu$ where $\tilde{\varphi}$ is a similar function that transitions in the interval $[3R_0,\,4R_0]$:
$$\tilde{\varphi}(r) = \xi\left(\frac{r-3R_0}{R_0}\right)\varphi_1 +  \left(1-\xi\left(\frac{r-3R_0}{R_0}\right)\right)\varphi_2.$$
We define
$${\bf W} = ({\bf U},\, \tilde{\bf U}).$$

\subsection{Construction of the Coefficients}
Construct $\tilde{f},\,\tilde{h}$ and $\tilde{\eta}$ in the exact same way as in Sections \ref{BuildingBlock} and \ref{Coupling}, for the function $\tilde{\varphi}$. We take 
$$ A_0 = 
\left(\begin{array}{cc|cc|cc|cc}
 f & 0 & 0 & \eta & 0 & 0 & 0 & 0 \\
 0 & h & -\eta & 0 & 0 & 0 & 0 & 0 \\
 \hline
 0 & -\eta & f & 0 & 0 & 0 & 0 & 0 \\
 \eta & 0 & 0 & h & 0 & 0 & 0 & 0 \\
 \hline
 0 & 0 & 0 & 0 & \tilde{f} & 0 & 0 & \tilde{\eta} \\
 0 & 0 & 0 & 0 & 0 & \tilde{h} & -\tilde{\eta} & 0 \\
 \hline
 0 & 0 & 0 & 0 & 0 & -\tilde{\eta} & \tilde{f} & 0 \\
 0 & 0 & 0 & 0 & \tilde{\eta} & 0 & 0 & \tilde{h}
\end{array}\right)(r)
$$
with respect to the coordinate system
$$(\nu,0,0,0),\, (\tau,0,0,0),\, (0,\nu,0,0),\, (0,\tau,0,0),...,(0,0,0,\tau),$$
where $(v,w,x,y)$ denotes the $4 \times 2$ matrix with rows $v,\,w,\,x$ and $y$.
Then $A_0$ is smooth and uniformly elliptic. (Indeed, the top left and lower right blocks are uniformly elliptic by the computations in Section \ref{Coupling}). Furthermore, we have
$$\text{div}(A_0(x)D{\bf W}) = \frac{1}{2}D{\bf W} \cdot x.$$

\subsection{Showing the Coefficients Depend Smoothly on ${\bf W}$}
We show that $A_0(x)$ can be written as $A({\bf W})$ for a uniformly elliptic, smooth matrix field $A$ on $\mathbb{R}^4$.

\begin{figure}
 \centering
    \includegraphics[scale=0.35]{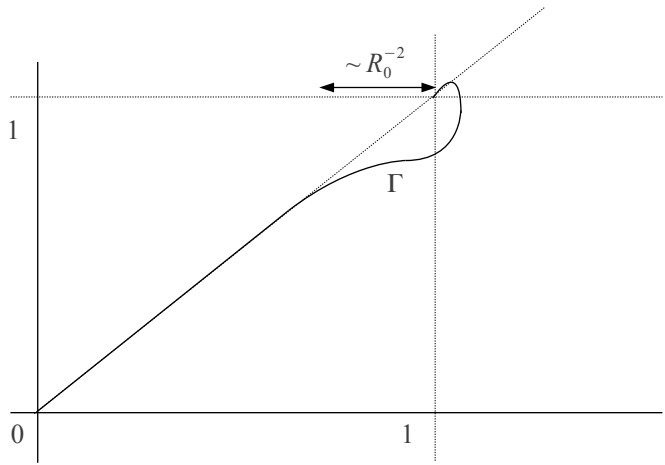}
 \caption{The image of $(\varphi,\, \tilde{\varphi})$ is a smooth embedded curve $\Gamma$.}
 \label{Gamma}
\end{figure}

Let $\Gamma \subset \mathbb{R}^2$ be the image $(\varphi,\tilde{\varphi})((0,\infty))$. Then $\Gamma$ is a smooth embedded curve consisting of two
segments on the diagonal $\theta = \frac{\pi}{4}$ connected by a short piece below the diagonal (see Figure \ref{Gamma}).

Define smooth functions $N$ and $H$ on $\Gamma$ by
$$N(\varphi(r), \,\tilde{\varphi}(r)) = \eta(r), \quad H(\varphi(r),\,\tilde{\varphi}(r)) = h(r).$$
Also, let
$$F(\varphi(r)) = f(r)$$
be a function on $[0, \max \varphi]$.
This definition makes sense because $f(r)$ is constant where $\varphi(r) \geq 1- \delta$ for some small $\delta$ (after possibly making $f$ transition to constant faster near $r = R_0$). One can extend
$F$ to a smooth, positive, bounded, even function $\mathcal{F}$ on $\mathbb{R}$ by letting $\mathcal{F}(s) = f(R_0)$ for $s \geq 1$, and by noticing that the expansion of $f$ near the origin has only even powers.

By construction we have that $N = 0$ on $\Gamma$ except for in a small square $Q_{\bar{\delta}}(1,1)$ of side length $2\bar{\delta}$ centered at $(1,1)$  (here $\bar{\delta}$ is of order $R_0^{-2}$).
Furthermore, $\mathcal{F}(x)$ is of order $R_0^2\log R_0$ for $(x,y) \in Q_{\bar{\delta}(1,1)}$. Note that $N$ is constant very close to $(1,1)$ on $\Gamma$. Extend $N$ to a smooth function $\mathcal{N}(x,y)$ on the positive quadrant that 
is less than order $\log R_0$ in $Q_{\bar{\delta}}(1,1)$ and vanishes outside of $Q_{\bar{\delta}(1,1)}$.

Next, we observe that $H = 1/2$ on $\Gamma$ away from $Q_{\bar{\delta}}$, and that near $(1,1)$ we have by construction that $H$ agrees with the function $4f(R_0) - \frac{4 f(R_0) - 1/2}{x}$. Extend $H$ to a smooth
function $\mathcal{H}$ on the positive quadrant that is identically $1/2$ away from $Q_{\bar{\delta}}$, and at least $1/3$ in the square.

For $(p,\,q) \in \mathbb{R}^4$, the functions $\mathcal{F}(|p|), \, \mathcal{H}(|p|,\,|q|)$ and $\mathcal{N}(|p|,\,|q|)$ are smooth. Define
$$A_{12}(p,\,q) = -A_{21}(p,\,q) = \mathcal{N}(|p|,\,|q|) \left(\begin{array}{cc}
 0 & 1 \\
 -1 & 0
 \end{array}\right),$$
 and
 \begin{align*}
 A_{11}(p,\,q) = A_{22}(p,\,q) &= \mathcal{F}(|p|) \frac{p \otimes p}{|p|^2} + \mathcal{H}(|p|,|q|) \frac{p^{\perp} \otimes p^{\perp}}{|p|^2} \\
  &= \frac{1}{2} I + (\mathcal{F}-1/2)(|p|)\frac{p \otimes p}{|p|^2} + (\mathcal{H}-1/2)(|p|,|q|) \frac{p^{\perp} \otimes p^{\perp}}{|p|^2}.
 \end{align*}
Then $A^{ij}_{\alpha\beta}|_{\alpha,\,\beta \leq 2}$ is a smooth, bounded, uniformly elliptic matrix field on $\mathbb{R}^4$. Indeed, $\mathcal{H} - 1/2$ is zero except for $(|p|,\,|q|)$ near $(1,1)$ and is larger than $-1/6$, and 
 $\mathcal{F} - 1/2$ is a smooth positive bounded function that vanishes on $\{p = 0\}$ and is of order $R_0^2 \log R_0$ where $\mathcal{N}$ is of order $\log R_0$.
 
Finally, it is clear from the definitions of $\mathcal{F},\, \mathcal{H}$ and $\mathcal{N}$ that $A^{ij}_{\alpha\beta}({\bf W}(x))|_{\alpha,\,\beta \leq 2}$ agree with the same components of $A_0(x)$.
 
Using a very similar procedure with $\tilde{f},\,\tilde{h}$ and $\tilde{\eta}$, one can also define uniformly elliptic smooth coefficients $A^{ij}_{\alpha\beta}|_{\alpha,\,\beta \geq 3}$ on $\mathbb{R}^4$ so that
$A^{ij}_{\alpha\beta}({\bf W}(x))|_{\alpha,\,\beta \geq 3}$ agree with the same components of $A_0(x)$. 
Taking the remaining coefficients to be zero completes the construction.

\begin{proof}[{\bf Proof of Theorem \ref{nonlinearmain}}]
We have constructed a smooth bounded map ${\bf W} : \mathbb{R}^2 \rightarrow \mathbb{R}^4$ and smooth uniformly elliptic coefficients $A$ on $\mathbb{R}^4$ verifying the hypotheses in Remark \ref{QuasilinearReduction} and 
the estimates (\ref{BoundaryConditions}).
\end{proof}

\section{Liouville Theorem}\label{Appendix}
In the final section we prove a Liouville theorem showing why $|{\bf U}|$ can not be radially increasing in two dimensions.

\begin{thm}\label{Liouville}
Any global, bounded solution ${\bf U} : \mathbb{R}^2 \rightarrow \mathbb{R}^m$ to the uniformly elliptic system
$$\text{div}(A(x)D{\bf U}) = f(x)D{\bf U} \cdot x$$
such that $f \geq 0$ and $|{\bf U}|$ is radially increasing is constant.
\end{thm}

\begin{rem}
The examples of Giusti-Miranda \cite{GM} and Star\`{a}-John \cite{SJ} show that the condition $n = 2$ is necessary. 
\end{rem}

\begin{proof}
The key observation is that, since $|{\bf U}|$ is radially increasing, we have
 $$0 \leq \frac{1}{2}f(x)\nabla |{\bf U}|^2 \cdot x = f(x) {\bf U} \cdot (D{\bf U} \cdot x).$$
In particular,
$$0 \leq \text{div}(AD{\bf U}) \cdot {\bf U}\psi^2$$
for any compactly supported $H^1$ function $\psi$. Integrating by parts and using uniform ellipticity one obtains the Caccioppoli inequality
$$\int_{\mathbb{R}^2} |D{\bf U}|^2\psi^2\,dx \leq C\int_{\mathbb{R}^2} |{\bf U}|^2|\nabla \psi|^2\,dx.$$
Since ${\bf U}$ is bounded we thus have
$$\int_{\mathbb{R}^2} |D{\bf U}|^2\psi^2\,dx \leq C \int_{\mathbb{R}^2} |\nabla \psi|^2\,dx.$$
Taking $\psi = 1$ in $B_1$, zero outside of $B_R$, and 
$$\psi = 1-\frac{\log r}{\log R} \quad \text{ for } 1 \leq r \leq R$$
the above inequality becomes
$$\int_{B_1} |D{\bf U}|^2\,dx \leq \frac{C}{\log R}.$$
Taking $R \rightarrow \infty$ we conclude that ${\bf U}$ is constant in $B_1$, and by a simple scaling argument that ${\bf U}$ is constant globally.
\end{proof}

\begin{rem}
By inspection of the proof, a Liouville theorem holds for any uniformly elliptic system in two dimensions of the form
$$\text{div}(A(x)D{\bf U}) = {\bf V} + g(|x|) D{\bf U} \cdot x^{\perp}$$ 
such that ${\bf V} \cdot {\bf U} \geq 0$. Indeed, after taking the dot product with ${\bf U}$, the last term becomes an angular derivative of $|{\bf U}|^2$,
which disappears when we multiply by a radially symmetric cutoff and integrate.

Such systems arise by searching for self-similar solutions to (\ref{ParabolicSystem}) with radially increasing modulus, that are invariant under rescalings
that e.g. fix $-\epsilon$-homogeneous maps (giving the term ${\bf V} = \frac{1}{2}(D{\bf U} \cdot x + \epsilon {\bf U})$) or have ``spiraling'' behavior (giving a 
term involving the angular derivative of ${\bf U}$).
\end{rem}





\section*{Acknowledgment}
This work was supported by NSF grant DMS-1501152. I thank A. Figalli and A. Vasseur for discussions.




\end{document}